\input ./style/arxiv-general.cfg
\documentclass[MSNbibl,number,citesort,dvips]{arxbj}
\makeatletter
   \@ifpackageloaded{graphicx}{}{\usepackage{graphicx}}
\makeatother

\volume{23}
\issue{1}
\pubyear{2017}
\firstpage{459}
\lastpage{478}
\doi{10.3150/15-BEJ749}
\docsubty{FLA}

\makeatletter
\def\sfrac#1#2{#1/#2}

\newcommand{\lleft}{\left}
\newcommand{\rrvert}{\vert}
\newcommand{\rright}{\right}
\newcommand{\rrVert}{\Vert}
\newcommand{\llvert}{\vert}
\newcommand{\llVert}{\Vert}
\newcommand{\eqref}[1]{(\ref{#1})}
\newcommand{\Ex}{\mathrm{E}}
\newcommand{\Var}{\operatorname{Var}}
\newcommand{\trace}{\operatorname{tr}}
\newcommand{\rank}{\operatorname{rank}}
\newtheorem{proposition}{Proposition}
\newtheorem{corollary}{Corollary}
\newtheorem{lemma}{Lemma}
\newremark{remark}{Remark}
\newcommand{\V}{{\mathsf{V}}}
\newcommand{\tsl}{T_{\lambda}}
\newcommand{\msl}{M_{\lambda}}
\newcommand{\qsl}{Q_{\lambda}}
\newcommand{\dmax}{d_{\max}}
\newcommand{\psimax}{\psi_{\max}}
\newcommand{\hsl}{H_{\lambda}}
\newcommand{\lmin}{\lambda_{\min}}
\makeatother

\begin{document}
\begin{frontmatter}

\title{Convergence analysis of block Gibbs samplers for Bayesian
linear mixed models with $p>N$}
\runtitle{Convergence of the Gibbs sampler}

\begin{aug}
\author{\inits{T.}\fnms{Tavis}~\snm{Abrahamsen}\corref{}\thanksref{e1}\ead[label=e1,mark]{tabrahamsen@stat.ufl.edu}}
\and
\author{\inits{J.P.}\fnms{James P.}~\snm{Hobert}\thanksref{e2}\ead[label=e2,mark]{jhobert@stat.ufl.edu}}
\address{Department of Statistics, University of Florida,
Gainesville, Florida 32611, USA.\\ \printead{e1,e2}}
\end{aug}

%
\received{\smonth{2} \syear{2015}}

%
\begin{abstract}
Exploration of the intractable posterior distributions associated
with Bayesian versions of the general linear mixed model is often
performed using Markov chain Monte Carlo. In particular, if a
conditionally conjugate prior is used, then there is a simple
two-block Gibbs sampler available. Rom{\'a}n and Hobert
[\textit{Linear Algebra Appl.} \textbf{473} (2015) 54--77] showed
that, when the priors are proper and the $X$ matrix has full column
rank, the Markov chains underlying these Gibbs samplers are nearly
always geometrically ergodic. In this paper, Rom{\'a}n and Hobert's
(2015) result is extended by allowing improper priors on the
variance components, and, more importantly, by removing all
assumptions on the $X$ matrix. So, not only is $X$ allowed to be
(column) rank deficient, which provides additional flexibility in
parameterizing the fixed effects, it is also allowed to have more
columns than rows, which is necessary in the increasingly important
situation where $p > N$. The full rank assumption on $X$ is at the
heart of Rom{\'a}n and Hobert's (2015) proof. Consequently, the
extension to unrestricted $X$ requires a substantially different
analysis.
\end{abstract}

%
\end{frontmatter}

\section{Introduction}
\label{sec:intro}

The general linear mixed model (GLMM) is one of the most frequently
applied statistical models. A GLMM with $r$ random factors takes the
form
\[
Y = X \beta+ \sum_{i=1}^r
Z_i u_i + e,
\]
where $Y$ is an observable $N \times1$ data vector, $X$ and
$\{Z_i\}_{i=1}^r$ are known matrices, $\beta$ is an unknown $p \times
1$ vector of regression coefficients, $\{u_i\}_{i=1}^r$ are
independent random vectors whose elements represent the various levels
of the random factors in the model, and $e \sim\mbox{N}_N(0,
\lambda^{-1}_e I)$. Assume that $e$ and $u:=  ( u_1^T~~u_2^T~~\cdots~~ u_r^T  )^T$ are independent, and that $u \sim\mbox{N}_q(0,
\Lambda^{-1})$, where $u_i$ is $q_i \times1$, $q = q_1 + \cdots+
q_r$, and $\Lambda= \bigoplus_{i=1}^r \lambda_{u_i} I_{q_i}$. Letting
$Z = (Z_1~~Z_2~~\cdots~~Z_r)$, we can write $\sum_{i=1}^r Z_i u_i =
Zu$. Let $\lambda$ denote the vector of precision parameters, i.e.,
$\lambda= (\lambda_e~~\lambda_{u_1}~~\cdots~~\lambda_{u_r})^T$. To rule out degenerate cases, we assume throughout
that $N \ge2$, and that $q_i \ge2$ for each $i=1,2,\dots,r$. For a
book-length treatment of the GLMM, which is sometimes called the
\textit{variance components model}, see~\cite{searcasemccu1992}.

In the Bayesian setting, prior distributions are assigned to the
unknown parameters, $\beta$ and $\lambda$. Unfortunately, the Bayes
estimators associated with any non-trivial prior cannot be obtained in
closed form. This is because such estimators take the form of ratios
of high-dimensional, intractable integrals. The dimensionality also
precludes the use of classical Monte Carlo methods that require the
ability to draw samples directly from the posterior distribution.
Instead, the parameter estimates are typically obtained using Markov
chain Monte Carlo (MCMC) methods. In particular, when (proper or
improper) conditionally conjugate priors are adopted for $\beta$ and
$\lambda$, there is a simple block Gibbs sampler that can be used to
explore the intractable posterior density. Let $\theta= (\beta^T~~u^T)^T$, and denote the posterior density as $\pi(\theta,\lambda|y)$,
where $y$ denotes the observed data vector. (Since $u$ is
unobservable, it is treated like a parameter.) When the conditionally
conjugate priors are adopted, it is straightforward to simulate from
$\theta|\lambda,y$, and from $\lambda|\theta,y$. Indeed,
$\theta|\lambda,y$ is multivariate normal and, given $(\theta,y)$, the
components of $\lambda$ are independent gamma variates. Hence, it is
straightforward to simulate a Markov chain,
$\{(\theta_n,\lambda_n)\}_{n=0}^\infty$, that has
$\pi(\theta,\lambda|y)$ as its invariant density. Our main results
concern the convergence properties of this block Gibbs sampler. We
now provide some background about Markov chains on $\mathbb{R}^d$,
which will allow us to describe our results and their practical
importance.

Let $V = \{V_m\}_{m=0}^\infty$ denote a Markov chain with state space
$\V\subset\mathbb{R}^d$ and assume the chain is Harris ergodic;
that is, $\psi$-irreducible, aperiodic and positive Harris recurrent (see
\cite{meyntwee1993} for definitions). Assume further that the
chain has a Markov transition density (with respect to Lebesgue
measure), $k: \V\times\V\rightarrow[0,\infty)$. Then, for any
measurable set $A$, we have
\[
\Pr(V_{m+1} \in A | V_m=v) = \int_A
k\bigl(v'|v\bigr) \,dv'.
\]
For $m \in\{2,3,4,\dots\}$, the $m$-step Markov transition density
(Mtd) is defined inductively as follows
\[
k^m\bigl(v'|v\bigr) = \int_\V
k^{m-1}\bigl(v'|u\bigr) k(u|v) \,du.
\]
Of course, $k^1 \equiv k$, and $k^m(\cdot|v)$ is the density of $V_m$
conditional on $V_0=v$. Suppose that the invariant probability
distribution also has a density (with respect to Lebesgue measure),
$\kappa: \V\rightarrow[0,\infty)$. The chain $V$ is
\textit{geometrically ergodic} if there exist $M: \V\rightarrow
[0,\infty)$ and $\gamma\in[0,1)$ such that, for all $m \in
\mathbb{N}$,
%
%
\begin{equation}
\label{eq:ge} \int_\V\bigl\llvert k^m\bigl(
v'| v\bigr) - \kappa\bigl(v'\bigr) \bigr\rrvert
\,dv' \le M(v) \gamma ^m.
\end{equation}
Of course, the quantity on the left-hand side of \eqref{eq:ge} is the
total variation distance between the invariant distribution and the
distribution of $V_m$ given $V_0=v$.

There are many important practical and theoretical benefits of using a
geometrically ergodic Markov chain as the basis of one's MCMC
algorithm (see, e.g., \cite{roberose1998,jonehobe2001,flegharajone2008}). Perhaps
the most important of these is the ability to construct valid
asymptotic standard errors for MCMC-based estimators. Let $h: \V
\rightarrow\mathbb{R}$ be a function such that $\int_\V|h(v)|
\kappa(v) \,dv < \infty$, and suppose that the chain $V$ is to serve
as the basis of an MCMC algorithm for estimating $\omega:= \int_\V
h(v) \kappa(v) \,dv < \infty$. Harris ergodicity guarantees that
the standard estimator of $\omega$, $\overline{h}_m:= \frac{1}{m}
\sum_{i=0}^{m-1} h(V_i)$, is\vspace*{1pt} strongly consistent. However, Harris
ergodicity is \textit{not} enough to ensure that $\overline{h}_m$
satisfies a central limit theorem (CLT). On the other hand, if $V$ is
geometrically ergodic and there exists an $\varepsilon>0$ such that
$\int_\V|h(v)|^{2+\varepsilon} \kappa(v) \,dv < \infty$, then
$\overline{h}_m$ does indeed satisfy a $\sqrt{m}$-CLT; that is, under
these conditions, there exists a positive, finite $\sigma^2$ such
that, as $m \rightarrow\infty$, $\sqrt{m} (\overline{h}_m -
\omega) \stackrel{d}{\rightarrow} \mbox{N}(0,\sigma^2)$. This is
extremely important from a practical standpoint because all of the
standard methods of calculating valid asymptotic standard errors for
$\overline{h}_m$ are based on the existence of this CLT
(Flegal, Haran and Jones \cite{flegharajone2008}).

There have been several studies of the convergence properties of the
block Gibbs sampler for the GLMM
(Johnson and Jones \cite{johnjone2010}, Rom{\'a}n and Hobert \cite{romahobe2012,romahobe2015}). These have
resulted in easily-checked sufficient conditions for geometric
ergodicity of the underlying Markov chain. However, in all of the
studies to date, the matrix $X$ has been assumed to have full column
rank. In this paper, we extend the results to the case where $X$ is
completely unrestricted. So, not only do we allow for a rank
deficient $X$, which provides additional flexibility in parameterizing
the fixed effects, we also allow for $X$ with $p > N$, which is
necessary in the increasingly important situation where there are more
predictors than data points. Two different families of conditionally
conjugate priors are considered, one proper and one improper. We now
describe our results, beginning with the results for proper priors.

Assume that $\beta$ and the components of $\lambda$ are all \textit{a
priori} independent, and that $\beta\sim
\mbox{N}_p(\mu_\beta,\Sigma_\beta)$, $\lambda_e \sim
\operatorname{Gamma}(a_0,b_0)$ and, for $i=1,\dots,r$, $\lambda_{u_i} \sim
\operatorname{Gamma}(a_i,b_i)$. Our result for proper priors, which is a
corollary of Proposition~\ref{prop:ge_p_pred} from
Section~\ref{sec:ana}, is as follows.

%
\begin{corollary}\label{cor:ge_p}\label{cor1}
Under a proper prior, the block Gibbs Markov chain,
$\{(\lambda_n,\theta_n)\}_{n=0}^\infty$, is geometrically ergodic if:
\begin{longlist}[2.]
\item[1.] $a_0 > \frac{1}{2}(\rank(Z)-N+2) $, and

\item[2.]$\min\{ a_1+\frac{q_1}{2}, \dots, a_r+\frac{q_r}{2}
\} > \frac{1}{2}  ( q-\rank(Z)  ) + 1 $.
\end{longlist}
\end{corollary}

The conditions of Corollary~\ref{cor:ge_p} are quite weak in
the sense that they would nearly always be satisfied in practice.
Indeed, it would typically be the case that $\rank(Z)-N < -2$ (making
the first condition vacuous) and $q-\rank(Z)$ is close to zero (making
the second condition easily satisfied). In fact, if $q=\rank(Z)$,
which is the case for many standard designs, then the second condition
is also vacuous.

Rom{\'a}n and Hobert \cite{romahobe2015} (hereafter R\&H15) proved this same result
under the restrictive assumption that $X$ has full column rank.
Moreover, the rank assumption is at the very heart of their proof.
Indeed, these authors established a geometric drift condition for the
marginal chain, $\{\theta_n\}_{n=0}^\infty$, but their drift
(Lyapunov) function is \textit{only valid when $X$ has full column
rank}. Our proof is significantly different. We analyze the other
marginal chain, $\{\lambda_n\}_{n=0}^\infty$, using a drift function
that does not involve the matrix $X$. Generally speaking, minor
changes in a drift function often lead to significant differences in
what one is able to prove. Thus, it is somewhat surprising that we
are able to recover \textit{exactly} the conditions of R\&H15. To be
fair, we are able to use several of their matrix bounds, but only
\textit{after} extending them to the case where $X$ is unrestricted.

When $X$ does not have full column rank, a flat prior on $\beta$ leads
to an improper posterior. Thus, the improper priors that we consider
are actually partially proper. In particular, assume again that
$\beta$ and the components of $\lambda$ are all \textit{a priori}
independent, and that $\beta\sim\mbox{N}_p(\mu_\beta,\Sigma_\beta)$.
But now take the prior on $\lambda_e$ to be (proportional to)
$\lambda_e^{a_0-1} e^{-b_0 \lambda_e}
I_{(0,\infty)}(\lambda_e)$, and for $i=1,\dots,r$, take the prior on
$\lambda_{u_i}$ to be $\lambda_{u_i}^{a_i-1} e^{-b_i \lambda_{u_i}}
I_{(0,\infty)}(\lambda_{u_i})$. Assume that $\min\{a_i,b_i\} \le0$
for at least one $i=0,1,\dots,r$; otherwise, we are back to the proper
priors described above. (See \cite{gelm2006} for a comprehensive
discussion about improper priors for variance components.) Let $W =
(X~~Z)$, so that $W \theta= X\beta+ Zu$, and define $\mathrm{SSE} =
\llVert  y - W \hat{\theta}\rrVert ^2$, where $\hat{\theta}
= (W^T W)^+ W^T y$
and superscript ``$+$'' on a matrix denotes Moore--Penrose inverse.
Our result for improper priors, which is another corollary of
Proposition~\ref{prop:ge_p_pred}, is as follows.

%
\begin{corollary}
\label{cor:ge_i}
Under an improper prior, the block Gibbs Markov chain,
$\{(\lambda_n,\theta_n)\}_{n=0}^\infty$, is geometrically ergodic if:
\begin{longlist}[2.]
\item[1.] $2 b_0 + \mathrm{SSE} > 0 $,

\item[2.] For each $i \in\{1,2,\dots,r\}$, either $b_i>0$ or $a_i<b_i=0$,

\item[3.]$a_0 > \frac{1}{2}(\rank(Z)-N+2) $, and

\item[4.]$\min\{ a_1+\frac{q_1}{2}, \dots, a_r+\frac{q_r}{2}
\} > \frac{1}{2}  ( q-\rank(Z)  ) + 1 $.
\end{longlist}
\end{corollary}

Note that the two conditions of Corollary~\ref{cor:ge_p} are
exactly the same as the third and fourth conditions of
Corollary~\ref{cor:ge_i}. Furthermore, the first two conditions of
Corollary~\ref{cor:ge_i} are \textit{necessary} for posterior
propriety \cite{suntsuthe2001}, and hence for geometric
ergodicity. Consequently, the commentary above regarding the weakness
of the conditions of Corollary~\ref{cor:ge_p} applies here as well.

Corollary~\ref{cor:ge_i} is the first convergence rate result for the
block Gibbs sampler for this set of partially proper priors.
Rom{\'a}n and Hobert \cite{romahobe2012} (hereafter R\&H12) proved a similar result (see
their Corollary 1) for a different family of improper priors in which
our proper multivariate normal prior on $\beta$ is replaced by a flat
prior. Of course, because they used a flat prior on $\beta$, their
results are only relevant in the case where $X$ has full column rank.

The remainder of this paper is organized as follows. A formal
definition of the block Gibbs Markov chain is given in
Section~\ref{sec:sampler}. Section~\ref{sec:ana} contains our
convergence rate analysis of the block Gibbs sampler under proper and
improper priors. A short discussion concerning an alternative result
for proper priors appears in Section~\ref{sec:dis}. Some technical
details are relegated to an \hyperref[append]{Appendix}.

\section{The block Gibbs sampler}
\label{sec:sampler}

The block Gibbs sampler is driven by the Markov chain
$\{(\theta_n,\lambda_n)\}_{n=0}^\infty$, which lives on the space
$\mathbb{R}^{p+q} \times\mathbb{R}_+^{r+1}$, where $\mathbb{R}_+:=
(0,\infty)$. The Markov transition density (of the version that
updates $\theta$ first) is given by
\[
k(\tilde{\theta}, \tilde{\lambda} | \theta, \lambda) = \pi(\tilde{\lambda}|
\tilde{\theta},y) \pi(\tilde{\theta }|\lambda,y).
\]
We will often suppress dependence on $y$, as we have in the Markov
transition density. The conditional densities,
$\pi(\lambda|\theta,y)$ and $\pi(\theta|\lambda,y)$, are now
described. The following formulas hold for both sets of priors
(proper and improper). The components of $\lambda$ are conditionally
independent given $\theta$, and we have
%
%
\begin{equation}
\label{eq:postlam_e} \lambda_{e} | \theta\sim\operatorname{Gamma}
\biggl(a_0 + \frac{N}{2}, b_0 + \frac{\llVert  y-W\theta\rrVert ^2}{2}
\biggr),
\end{equation}
and, for $i=1,\ldots,r$,
%
%
\begin{equation}
\label{eq:postlam_i} \lambda_{u_i} | \theta\sim\operatorname{Gamma}
\biggl(a_{i} + \frac
{q_{i}}{2}, b_{i} + \frac{\llVert  u_{i}\rrVert ^2}{2}
\biggr).
\end{equation}
When considering improper priors, we assume that these conditional
distributions are all well defined. In other words, we assume that
$\{a_i\}_{i=0}^r$ and $\{b_i\}_{i=0}^r$ are such that all of the shape
and rate parameters in the gamma distributions above are strictly
positive. Of course, this is not enough to guarantee posterior
propriety. However, the drift technique that we employ is equally
applicable to positive recurrent (proper posterior) and non-positive
recurrent (improper posterior) Markov chains \cite{roma2012}.
Furthermore, geometrically ergodic chains are necessarily positive
recurrent, so any Gibbs Markov chain that we conclude is geometrically
ergodic, necessarily corresponds to a proper posterior. Consequently,
there is no need to check for posterior propriety before proceeding
with the convergence analysis.\vspace*{1pt}

Now define $\tsl= \lambda_{e} X^{T}X + \Sigma_\beta^{-1}$, $\msl=
I - \lambda_{e}
X T_{\lambda}^{-1}X^{T}$, and $\qsl= \lambda_{e} Z^{T} \msl Z
+\Lambda$.
Conditional on $\lambda$, $\theta$ is multivariate normal with mean
%
%
\begin{equation}
\label{eq:postmean} \Ex[\theta| \lambda] = \lleft[ %
\begin{array} {c}
T_{\lambda}^{-1}\bigl(\lambda_{e} X^{T}y +
\Sigma_\beta^{-1}\mu_{\beta
}\bigr) -
\lambda_{e}^{2} T_{\lambda}^{-1}
X^{T} Z Q_{\lambda}^{-1}Z^{T} \bigl(\msl y -
X T_{\lambda}^{-1}\Sigma _\beta^{-1}
\mu_{\beta}\bigr)
\\[4pt]
\lambda_e Q_{\lambda}^{-1}Z^T \bigl(\msl
y - X T_{\lambda}^{-1}\Sigma _\beta^{-1}
\mu_{\beta}\bigr) \end{array} %
 \rright],
\end{equation}
and covariance matrix
%
%
\begin{equation}
\label{eq:postvar} \Var[\theta|\lambda] = \lleft[ %
\begin{array} {cc}
T_{\lambda}^{-1}+ \lambda_{e}^{2}
T_{\lambda}^{-1}X^{T}Z Q_{\lambda
}^{-1}Z^{T}X
T_{\lambda}^{-1}& - \lambda_{e} T_{\lambda}^{-1}X^{T}Z
Q_{\lambda}^{-1}
\\[4pt]
- \lambda_{e} Q_{\lambda
}^{-1}Z^{T} X
T_{\lambda}^{-1}& Q_{\lambda}^{-1} \end{array}
 \rright].
\end{equation}
(A derivation of these conditionals can be found in
\cite{abra2015}.)

The two marginal sequences, $\{\theta_n\}_{n=0}^\infty$ and
$\{\lambda_n\}_{n=0}^\infty$, are themselves Markov chains, and it is
easy to establish that (when the posterior is proper) all three chains
are Harris ergodic. Moreover, geometric ergodicity is a solidarity
property for these three chains, that is, either all three chains are
geometrically ergodic, or none of them is (see, e.g., \cite{roberose2001,diackharsalo2008,romahobepres2014}).
Again, in contrast with R\&H15, who analyzed the $\theta$-chain,
$\{\theta_n\}_{n=0}^\infty$, we establish our results by analyzing the
$\lambda$-chain, $\{\lambda_n\}_{n=0}^\infty$. The Mtd of the
$\lambda$-chain is given by
\[
k_l(\tilde{\lambda} | \lambda) = \int_{\mathbb{R}^{p+q}} \pi(
\tilde{\lambda}|\theta,y) \pi(\theta|\lambda,y) \,d\theta.
\]
R\&H12 also analyzed the $\lambda$-chain, and their analysis serves as
a road map for ours. In fact, we use the same drift function as
R\&H12.

\section{Convergence analysis of the block Gibbs sampler}
\label{sec:ana}

In order to state our main result, we require a couple of definitions.
For $i=1,2,\dots,r$, let $R_i$ be the $q_{i} \times q$ matrix defined
as $R_i = [0_{q_i \times q_1} \cdots0_{q_i \times q_{i-1}}~~I_{q_i
\times q_i }~~0_{q_i\times q_{i+1}}\cdots0_{q_ i \times q_r}]$.
Note that $u_i = R_i u$. Let $P_{Z^T Z}$ denote the orthogonal
projection onto the column space of $Z^T Z$. Finally, define
\[
\tilde{s} = \min\biggl\{ a_0 + \frac{N}{2}, a_1 +
\frac{q_1}{2}, \dots, a_r + \frac{q_r}{2} \biggr\}.
\]
The following result holds for both sets of priors (proper and
improper).

%
\begin{proposition}
\label{prop:ge_p_pred}
The block Gibbs sampler Markov chain,
$\{\theta_{n},\lambda_{n}\}_{n=0}^{\infty}$, is geometrically
ergodic if:
\begin{longlist}[4.]
\item[1.]$\tilde{s} > 0 $;
\item[2.]$2 b_0 + \mbox{SSE} > 0 $;
\item[3.] For each $i \in\{1,2,\dots,r\}$, either $b_i>0$ or $a_i<b_i=0
$; and
\item[4.] There exists $s \in(0,1] \cap(0,\tilde{s})$ such that
\end{longlist}
%
%
\begin{eqnarray}\label{eq:hypoth}
&& \max\Biggl\{ \frac{\Gamma(a_0 + \sfrac{N}{2} -s)}{\Gamma(a_0 +
\sfrac{N}{2})} \biggl( \frac{\rank(Z)}{2}
\biggr)^s,
\nonumber\\[-8pt]\\[-8pt]\nonumber
&&\quad \sum_{i=1}^{r}
\frac{\Gamma(a_{i} + \sfrac{q_{i}}{2} - s)}{\Gamma(
a_{i} + \sfrac{q_{i}}{2})} \biggl( \frac{\trace(R_{i}(I -
P_{Z^{T}Z})R_{i}^{T})}{2} \biggr)^s \Biggr\} < 1.
\end{eqnarray}
\end{proposition}

%
\begin{remark}
When the prior is proper, that is, when $a_i>0$ and $b_i>0$ for all
$i \in\{0,1,\dots,r\}$, the first three conditions are
automatically satisfied, and $\tilde{s}>1$. On the other hand, when
the prior is improper, these three conditions ensure that
$\pi(\lambda|\theta,y)$ is well defined.
\end{remark}

Before embarking on our proof of Proposition~\ref{prop:ge_p_pred}, we
quickly demonstrate that Corollaries~\ref{cor:ge_p} and \ref{cor:ge_i}
follow immediately from it.

\begin{pf*}{Proof of Corollary~\ref{cor:ge_p}}
Since the prior is proper, it is enough to show that the conditions
of Corollary~\ref{cor:ge_p} imply that \eqref{eq:hypoth} is
satisfied for some $s \in(0,1]$. We show that this is indeed the
case, with $s=1$. First,
\[
\frac{\Gamma(a_0 + \sfrac{N}{2} -1)}{\Gamma(a_0 +
\sfrac{N}{2})} \frac{\rank(Z)}{2} = \frac{\rank(Z)}{2a_0
+ N - 2},
\]
which yields the first half of \eqref{eq:hypoth}. Now note that
\[
\sum_{i=1}^{r} \trace\bigl(R_{i}(I
-P_{Z^{T}Z})R_{i}^{T}\bigr) = \trace \Biggl[ (I -
P_{Z^{T}Z}) \sum_{i=1}^{r}
R_{i}^{T} R_{i}) \Biggr] = \trace(I -
P_{Z^{T}Z}) = q - \rank(Z).
\]
Thus,
\begin{eqnarray*}
\sum_{i=1}^{r} \frac{\Gamma(a_{i} + \sfrac{q_{i}}{2} - 1)}{\Gamma(
a_{i} + \sfrac{q_{i}}{2})}
\frac{\trace(R_{i}(I -
P_{Z^{T}Z})R_{i}^{T})}{2} & =& \sum_{i=1}^{r}
\frac{\trace(R_{i}(I
- P_{Z^{T}Z})R_{i}^{T})}{2a_i + q_i - 2}
\\[-4pt]
& \le&\frac{q -
\rank(Z)}{\min_{i=1,\ldots,r} (2a_i+q_i-2)},
\end{eqnarray*}
which yields the second half of \eqref{eq:hypoth}.
\end{pf*}

\begin{pf*}{Proof of Corollary~\ref{cor:ge_i}}
First note that conditions 3 and 4 of Corollary~\ref{cor:ge_i}
imply that $\tilde{s}>1$. The rest of the proof is the same as the
proof of Corollary~\ref{cor:ge_p}.
\end{pf*}

Our proof of Proposition~\ref{prop:ge_p_pred} is based on four lemmas,
which are proven in the \hyperref[append]{Appendix}. Let $\dmax$ denote the largest
singular value of the matrix $\tilde{X} = X \Sigma_{\beta}^{1/2}$.

%
\begin{lemma}
\label{lem:ineq1}
For each $i \in\{1,2,\dots,r\}$, we have
\[
\trace\bigl(R_{i} Q_{\lambda}^{-1}R_{i}^{T}
\bigr) \leq\bigl(\dmax^2 + \lambda_{e}^{-1}\bigr)
\trace \bigl( R_{i} \bigl(Z^{T}Z\bigr)^+
R_{i}^{T} \bigr) + \trace \bigl( R_{i}(I -
P_{Z^{T}Z}) R_{i}^{T} \bigr) \sum
_{j=1}^{r} \lambda_{u_j}^{-1}.
\]
\end{lemma}

%
\begin{lemma}
\label{lem:ineq2}
$\trace( W \Var(\theta|\lambda) W^{T}) \leq\lambda_{e}^{-1}
\rank(Z) + \dmax^{2} \rank(Z) + \trace(X \Sigma_{\beta} X^T)$.
\end{lemma}

%
\begin{lemma}
\label{lem:ineq3}
There exist finite constants $K_1$ and $K_2$, not depending on
$\lambda$, such that $\llVert \Ex[R_i u | \lambda]\rrVert
\le\sqrt{q_i} K_1
$ for $i=1,\dots,r$, and $\llVert  y - W \Ex[\theta| \lambda
]\rrVert  \le
K_2$.
\end{lemma}

%
\begin{remark}
\label{rem:constants}
The constants $K_1$ and $K_2$ are defined in the \hyperref[append]{Appendix}. They do
not have a closed form.
\end{remark}

We will write $A \preceq B$ to mean that $B-A$ is nonnegative
definite. Let $\psimax$ denote the largest eigenvalue of $Z^T Z$.

%
\begin{lemma}
\label{lem:ineq4}
For each $i \in\{1,2,\dots,r\}$, we have $(\psimax\lambda_e +
\lambda_{u_i})^{-1} I_{q_i} \preceq R_i Q_{\lambda}^{-1}R_i^T $.
\end{lemma}

\begin{pf*}{Proof of Proposition~\ref{prop:ge_p_pred}}
Define the drift function as follows
\[
v(\lambda) = \alpha\lambda_{e}^{c} + \alpha
\lambda_{e}^{-s} + \sum_{i=1}^{r}
\lambda_{u_i}^{c} + \sum_{i=1}^{r}
\lambda_{u_i}^{-s},
\]
where $\alpha$ and $c$ are positive constants (that are explicitly
constructed in the proof), and $s$ is from the fourth condition in
Proposition~\ref{prop:ge_p_pred}. We will show that there exist $\rho
\in[0,1)$ and a finite constant $L$ such that
%
%
\begin{equation}
\label{eq:drift_cond} \Ex\bigl[v(\tilde{\lambda})|\lambda\bigr] = \int
_{\mathbb{R}_+^{r+1}} v(\tilde{\lambda}) k_l(\tilde{\lambda} |
\lambda) \,d\tilde{\lambda} \leq\rho v(\lambda) + L.
\end{equation}
Then because the $\lambda$-chain is a Feller chain \cite{abra2015}
and the function $v(\cdot)$ is unbounded off compact sets (R\&H12), by
Meyn and Tweedie's \cite{meyntwee1993}, Lemma 15.2.8, the geometric drift condition
\eqref{eq:drift_cond} implies that the $\lambda$-chain is
geometrically ergodic. We now establish \eqref{eq:drift_cond}.\vadjust{\goodbreak}

First, note that
\begin{eqnarray*}
\Ex\bigl[v(\tilde{\lambda})|\lambda\bigr] & =& \int_{\mathbb{R}^{p+q}}
\biggl[ \int_{\mathbb{R}_+^{r+1}} v(\tilde{\lambda}) \pi(\tilde{\lambda}|
\theta,y) \,d\tilde{\lambda} \biggr] \pi(\theta|\lambda,y) \,d\theta= \Ex \bigl[ \Ex
\bigl[ v(\tilde{\lambda})|\theta \bigr] | \lambda \bigr]
\\
& =& \alpha \Ex \bigl[ \Ex \bigl[ \tilde{\lambda}_e^c|
\theta \bigr] | \lambda \bigr] + \alpha\Ex \bigl[ \Ex \bigl[ \tilde{
\lambda}_e^{-s}|\theta \bigr] | \lambda \bigr]
\\
&&{} + \sum_{i=1}^r \Ex \bigl[ \Ex \bigl[
\tilde{\lambda}_{u_i}^c|\theta \bigr] | \lambda \bigr] + \sum
_{i=1}^r \Ex \bigl[ \Ex \bigl[ \tilde{
\lambda}_{u_i}^{-s}|\theta \bigr] | \lambda \bigr].
\end{eqnarray*}
Using \eqref{eq:postlam_e} and the fact that $0 < b_0+\mbox{SSE}/2 <
b_0 + \llVert  y-W\theta\rrVert ^2/2$, we have
%
%
\begin{equation}
\label{eq:term1} %
\begin{aligned} \Ex \bigl[ \tilde{
\lambda}_e^c|\theta \bigr] &=& \frac{\Gamma(a_0 +
\sfrac{N}{2} + c)} {\Gamma(a_0 + \sfrac{N}{2})} \biggl(
b_0 + \frac{\llVert  y-W\theta\rrVert ^2}{2} \biggr)^{-c}
\\
&\le& \frac{\Gamma(a_0 +
\sfrac{N}{2} + c)} {\Gamma(a_0 + \sfrac{N}{2})} \biggl( b_0 + \frac{\mbox{SSE}}{2}
\biggr)^{-c}. \end{aligned} %
\end{equation}
As we shall see, since this upper bound does not depend on $\theta$,
it can be absorbed into the constant term, $L$, and we will no longer
have to deal with this piece of the drift function. Now,
%
%
\begin{equation}
\label{eq:term2} %
\begin{aligned} \Ex \bigl[ \tilde{
\lambda}_e^{-s}|\theta \bigr] &=& \frac{\Gamma(a_0 +
\sfrac{N}{2} - s)} {\Gamma(a_0 + \sfrac{N}{2})} \biggl(
b_0 + \frac{\llVert  y-W\theta\rrVert ^2}{2} \biggr)^s
\\
&\le& \frac{\Gamma(a_0 +
\sfrac{N}{2} - s)} {\Gamma(a_0 + \sfrac{N}{2})} \biggl( |b_0|^s + \biggl[
\frac{\llVert  y-W\theta\rrVert ^2}{2} \biggr]^s \biggr), \end{aligned} %
\end{equation}
where the inequality follows from the fact that $(x_1+x_2)^{\xi} \leq
x_1^\xi+x_2^\xi$ for $x_1,x_2 \geq0$ and $\xi\in(0,1]$. Similarly,
using \eqref{eq:postlam_i}, for each $i \in\{1,\dots,r\}$ we have
%
%
\begin{equation}
\label{eq:term3} %
\begin{aligned} \Ex \bigl[ \tilde{
\lambda}_{u_i}^{-s}|\theta \bigr] &=& \frac{\Gamma
(a_i +
\sfrac{q_i}{2} - s)} {\Gamma(a_i + \sfrac{q_i}{2})} \biggl(
b_i + \frac{\llVert  u_{i}\rrVert ^2}{2} \biggr)^s
\\
&\le&\frac
{\Gamma(a_i + \sfrac{q_i}{2}
- s)} {\Gamma(a_i + \sfrac{q_i}{2})} \biggl( b_i^s + \biggl[
\frac{\llVert  u_{i}\rrVert ^2}{2} \biggr]^s \biggr). \end{aligned} %
\end{equation}
Now, for each $i \in\{1,\dots,r\}$, we have
%
%
\begin{equation}
\label{eq:term4} %
\begin{aligned} \Ex \bigl[ \tilde{
\lambda}_{u_i}^c|\theta \bigr] &=& \frac{\Gamma
(a_i +
\sfrac{q_i}{2} + c)} {\Gamma(a_i + \sfrac{q_i}{2})} \biggl(
b_i + \frac{\llVert  u_{i}\rrVert ^2}{2} \biggr)^{-c}
\\
&\le&\frac
{\Gamma(a_i +
\sfrac{q_i}{2} + c)} {\Gamma(a_i + \sfrac{q_i}{2})} \biggl[ \biggl( \frac{\llVert  u_{i}\rrVert ^2}{2} \biggr)^{-c}
I_{\{0\}}(b_i) + b_i^{-c}
I_{\mathbb{R}_+}(b_i) \biggr]. \end{aligned} %
\end{equation}
Note that when $b_i>0$ there is a simple upper bound for this term
that does not depend on $\theta$. Therefore, we will first consider
the case in which $\min\{b_1,\dots,b_r\}>0$, and we will return to
the other (more complicated) case later.

Assume for the time being that $\min\{b_1,\dots,b_r\}>0$. Then
combining \eqref{eq:term1}, \eqref{eq:term2}, \eqref{eq:term3} and
\eqref{eq:term4}, and applying Jensen's inequality twice yields
%
%
\begin{eqnarray}
\label{eq:main} \Ex\bigl[v(\tilde{\lambda})|\lambda\bigr] &\le&
\frac{\alpha}{2^s}
\frac{\Gamma(a_0 + \sfrac{N}{2} - s)} {\Gamma(a_0 + \sfrac{N}{2})}  \Ex \bigl[ \llVert y - W\theta\rrVert ^2 |
\lambda \bigr]^s
\nonumber\\[-8pt]\\[-8pt]\nonumber
&&{} + 2^{-s} \sum_{i=1}^r
\frac{\Gamma(a_i + \sfrac{q_i}{2} - s)}{
\Gamma(a_i + \sfrac{q_i}{2})} \Ex \bigl[ \llVert u_i\rrVert ^2 |
\lambda \bigr]^s + K_0,
\end{eqnarray}
where
\begin{eqnarray*}
K_0 &=& \alpha\frac{\Gamma(a_0 + \sfrac{N}{2} + c)} {\Gamma(a_0 +
\sfrac{N}{2})} \biggl( b_0 +
\frac{\mbox{SSE}}{2} \biggr)^{-c} + \alpha\frac{\Gamma(a_0 + \sfrac{N}{2} - s)} {\Gamma(a_0 +
\sfrac{N}{2})}
|b_0|^s
\\
&&{} + \sum_{i=1}^r \biggl[
\frac{\Gamma(a_i
+ \sfrac{q_i}{2} + c)} {\Gamma(a_i + \sfrac{q_i}{2})} b_i^{-c} + \frac{\Gamma(a_i + \sfrac{q_i}{2} - s)} {\Gamma(a_i + \sfrac{q_i}{2})}
b_i^s \biggr].
\end{eqnarray*}
It follows from \eqref{eq:postvar} that
\[
\Ex \bigl[ \llVert y - W\theta\rrVert ^2 | \lambda \bigr] = \trace
\bigl( W\Var(\theta|\lambda)W^T \bigr) + \bigl\llVert y-W\Ex[\theta|
\lambda]\bigr\rrVert ^{2}.
\]
Similarly, since $u_i = R_i u$, we also have
\[
\Ex \bigl[ \llVert u_i\rrVert ^2 | \lambda \bigr] = \Ex
\bigl[ \llVert R_i u\rrVert ^2 | \lambda \bigr] = \trace
\bigl(R_i Q_{\lambda}^{-1}R_i^{T}
\bigr) + \bigl\llVert \Ex [R_i u |\lambda]\bigr\rrVert ^2.
\]
Now, using Lemmas~\ref{lem:ineq2} and \ref{lem:ineq3}, we have
%
%
\begin{eqnarray}
\label{eq:t1} \Ex \bigl[ \llVert y - W\theta\rrVert ^2 | \lambda
\bigr]^s & \le& \bigl[ \lambda_{e}^{-1} \rank(Z) +
\dmax^{2} \rank(Z) + \trace\bigl(X \Sigma_{\beta}
X^{T}\bigr) + K_2^2 \bigr]^s
\nonumber\\[-8pt]\\[-8pt]\nonumber
& \le& \lambda_{e}^{-s} \bigl(\rank(Z)\bigr)^s +
\bigl[\dmax^{2} \rank(Z) + \trace\bigl(X \Sigma_{\beta}
X^{T}\bigr) + K_2^2 \bigr]^s.
\end{eqnarray}
Similarly, using Lemmas~\ref{lem:ineq1} and \ref{lem:ineq3}, we have
%
%
\begin{eqnarray}
\label{eq:t2} \Ex \bigl[ \llVert u_i\rrVert ^2 |
\lambda \bigr]^s & \le& \Biggl[ \bigl(\dmax^2 +
\lambda_{e}^{-1}\bigr) \trace \bigl( R_{i}
\bigl(Z^{T}Z\bigr)^+ R_{i}^{T} \bigr)
\nonumber
\\
&&{}+\trace \bigl( R_{i}(I - P_{Z^{T}Z})
R_{i}^{T} \bigr) \sum_{j=1}^{r}
\lambda_{u_j}^{-1} + q_i K_1^2
\Biggr]^s
\nonumber\\[-8pt]\\[-8pt]\nonumber
& \le &\lambda_{e}^{-s} \bigl(\trace \bigl( R_{i}
\bigl(Z^{T}Z\bigr)^+ R_{i}^{T} \bigr)
\bigr)^s + \sum_{j=1}^{r}
\lambda_{u_j}^{-s} \bigl( \trace \bigl( R_{i}(I -
P_{Z^{T}Z}) R_{i}^{T} \bigr) \bigr)^s
\nonumber
\\
&&{} + \bigl[ \dmax^2 \trace \bigl( R_{i}
\bigl(Z^{T}Z\bigr)^+ R_{i}^{T} \bigr) +
q_i K_1^2 \bigr]^s.\nonumber
\end{eqnarray}
Define a function $\delta(\cdot)$ as follows:
\[
\delta(\alpha) = \frac{\Gamma(a_0 + \sfrac{N}{2} - s)} {\Gamma(a_0 +
\sfrac{N}{2})} \biggl( \frac{\rank(Z)}{2}
\biggr)^s + \frac{1}{\alpha} \sum_{i=1}^r
\frac{\Gamma(a_i + \sfrac{q_i}{2} - s)} {\Gamma(a_i +
\sfrac{q_i}{2})} \biggl( \frac{\trace ( R_{i} (Z^{T}Z)^+ R_{i}^{T}
 )}{2} \biggr)^s.
\]
Combining \eqref{eq:main}, \eqref{eq:t1} and \eqref{eq:t2} yields
%
%
\begin{eqnarray}
\label{eq:main2} \Ex\bigl[v(\tilde{\lambda})|\lambda\bigr] &\le& \alpha\delta(\alpha)
\lambda_e^{-s}
\nonumber\\[-8pt]\\[-8pt]\nonumber
&&{} + \Biggl[ \sum
_{i=1}^r \frac{\Gamma(a_i +
\sfrac{q_i}{2} - s)} {\Gamma(a_i + \sfrac{q_i}{2})} \biggl(
\frac{\trace ( R_{i}(I - P_{Z^{T}Z}) R_{i}^{T}  )}{2} \biggr)^s \Biggr] \sum
_{j=1}^{r} \lambda_{u_j}^{-s} + L,
\end{eqnarray}
where
\begin{eqnarray*}
L &=& K_0 + \frac{\alpha}{2^s} \frac{\Gamma(a_0 + \sfrac{N}{2} - s)}{
\Gamma(a_0 + \sfrac{N}{2})}  \bigl[
\dmax^{2} \rank(Z) + \trace\bigl(X \Sigma_{\beta}
X^{T}\bigr) + K_2^2 \bigr]^s
\\
&&{} + 2^{-s} \sum_{i=1}^r
\frac{\Gamma(a_i + \sfrac{q_i}{2} - s)} {\Gamma(a_i + \sfrac{q_i}{2})} \bigl[ \dmax^2 \trace \bigl( R_{i}
\bigl(Z^{T}Z\bigr)^+ R_{i}^{T} \bigr) +
q_i K_1^2 \bigr]^s.
\end{eqnarray*}
Next, defining
\[
\rho(\alpha):= \max \Biggl\{ \delta(\alpha), \sum_{i=1}^r
\frac{\Gamma(a_i + \sfrac{q_i}{2} - s)}{
\Gamma(a_i + \sfrac{q_i}{2})} \biggl( \frac{\trace ( R_{i}(I -
P_{Z^{T}Z}) R_{i}^{T}  )}{2} \biggr)^s \Biggr\},
\]
we have from \eqref{eq:main2} that
\begin{eqnarray*}
\Ex\bigl[v(\tilde{\lambda})|\lambda\bigr] & \le&\alpha\rho(\alpha)
\lambda_e^{-s} + \rho(\alpha) \sum
_{j=1}^{r} \lambda_{u_j}^{-s} + L
\nonumber
\\
& \le&\rho(\alpha) \Biggl( \alpha\lambda_e^c + \alpha
\lambda_e^{-s} + \sum_{j=1}^{r}
\lambda_{u_j}^c + \sum_{j=1}^{r}
\lambda_{u_j}^{-s} \Biggr) + L
\nonumber
\\
& = &\rho(\alpha) v(\lambda) + L.
\end{eqnarray*}
Hence, all that is left is to demonstrate the existence of an $\alpha
\in(0,\infty)$ such that $\rho(\alpha) \in[0,1)$. By~\eqref{eq:hypoth}, we know that
\[
\sum_{i=1}^{r} \frac{\Gamma(a_{i} + \sfrac{q_{i}}{2} - s)}{\Gamma(
a_{i} + \sfrac{q_{i}}{2})} \biggl(
\frac{\trace(R_{i}(I -
P_{Z^{T}Z})R_{i}^{T})}{2} \biggr)^s < 1.
\]
Therefore, it suffices to show that there exists an $\alpha\in
(0,\infty)$ such that $\delta(\alpha) < 1$. But $\delta(\alpha) < 1$
as long as
%
%
\begin{equation}
\label{eq:alpha_bound} \alpha> \frac{\sum_{i=1}^{r} \frac{\Gamma(a_{i} + \sfrac{q_{i}}{2} -
s)}{\Gamma( a_{i} + \sfrac{q_{i}}{2})}  ( \sfrac{\trace(R_{i}(I
- P_{Z^{T}Z})R_{i}^{T})}{2}  )^s}{1-\frac{\Gamma(a_0 +
\sfrac{N}{2} -s)}{\Gamma(a_0 + \sfrac{N}{2})}  (
\sfrac{\rank(Z)}{2}  )^s},
\end{equation}
which is a well-defined positive number by \eqref{eq:hypoth}. The
result has now been proven for the case in which $\min
\{b_1,\dots,b_r\}>0$.

%
\begin{remark}
Note that the two terms in the drift function involving $c$ were
both absorbed into the constant in the first step of the iterated
expectation. It follows that, at least in the case where $\min
\{b_1,\dots,b_r\}>0$, any $c>0$ can be used in the drift function.
\end{remark}

We now proceed to the case in which there is at least one $b_i=0$.
Let $B = \{ i \in\{1,\dots,r\}: b_i = 0 \}$. It follows
from the development above that the following holds for any $c>0$:
%
%
\begin{equation}
\label{eq:for_later} \alpha\Ex \bigl[ \tilde{\lambda}_e^c|
\lambda \bigr] + \alpha\Ex \bigl[ \tilde{\lambda}_e^{-s}|
\lambda \bigr] + \sum_{i \notin B} \Ex \bigl[ \tilde{
\lambda}_{u_i}^c|\lambda \bigr] + \sum
_{i=1}^r \Ex \bigl[ \tilde{\lambda}_{u_i}^{-s}|
\lambda \bigr] \le\rho(\alpha) \Biggl( \alpha \lambda_e^{-s}
+ \sum_{j=1}^{r}\lambda_{u_j}^{-s}
\Biggr) + L.
\end{equation}
Of course, if $\alpha$ satisfies \eqref{eq:alpha_bound}, then
$\rho(\alpha) \in[0,1)$. Now suppose we can find $c>0$, $\alpha$
satisfying \eqref{eq:alpha_bound}, and $\rho'(\alpha) \in[0,1)$ such
that
%
%
\begin{equation}
\label{eq:suffice} \sum_{i \in B} \Ex \bigl[ \tilde{
\lambda}_{u_i}^c|\lambda \bigr] \le \rho'(
\alpha) \biggl( \alpha\lambda_e^c + \sum
_{i \in B} \lambda_{u_i}^c \biggr).
\end{equation}
Then combining \eqref{eq:for_later} and \eqref{eq:suffice}, we would
have
\begin{eqnarray*}
\Ex\bigl[v(\tilde{\lambda})|\lambda\bigr] & \le&\rho(\alpha) \Biggl( \alpha
\lambda_e^{-s} + \sum_{j=1}^{r}
\lambda_{u_j}^{-s} \Biggr) + L + \rho'(\alpha)
\biggl( \alpha\lambda_e^c + \sum
_{i \in B} \lambda_{u_i}^c \biggr)
\\[-2pt]
& \le&\max\bigl\{ \rho(\alpha), \rho'(\alpha) \bigr\} v(\lambda) +
L,
\end{eqnarray*}
which establishes the drift condition. Therefore, to prove the result
when $\min\{b_1,\dots,b_r\} = 0$, it suffices to establish
\eqref{eq:suffice}. If $i \in B$, then
\[
\Ex \bigl[ \tilde{\lambda}_i^c|\theta \bigr] =
\frac{\Gamma(a_i +
\sfrac{q_i}{2} + c)} {\Gamma(a_i + \sfrac{q_i}{2})} \biggl( \frac{\llVert  u_{i}\rrVert ^2}{2} \biggr)^{-c}.
\]
It follows from \eqref{eq:postvar} that the conditional distribution
of $(R_i Q_{\lambda}^{-1}R_i^T)^{-1/2} u_i$ given $\lambda$ is multivariate
normal with identity covariance matrix. Thus, $u_i^T (R_i Q_{\lambda}^{-1}
R_i^T)^{-1} u_i$ has a non-central chi-squared distribution with $q_i$
degrees of freedom. An application of Lemma 4 from R\&H12 shows that,
if $c \in(0,1/2)$, then
\[
\Ex \bigl[ \bigl( u_i^T \bigl(R_i
Q_{\lambda}^{-1}R_i^T
\bigr)^{-1} u_i \bigr)^{-c} | \lambda \bigr]
\le2^{-c} \frac{\Gamma(\sfrac{q_i}{2} -
c)}{\Gamma(\sfrac{q_{i}}{2})}.
\]
Putting this together with Lemma~\ref{lem:ineq4}, we have that, if $i
\in B$ and $c \in(0,1/2)$, then
\begin{eqnarray*}
\Ex \bigl[\bigl(\|u_i\|^{2}\bigr)^{-c} |
\lambda \bigr] & =& (\psimax\lambda_e + \lambda_{u_i})^c
\Ex \bigl[ \bigl( u_i^T (\psimax\lambda_e +
\lambda_{u_i}) I_{q_i} u_i \bigr)^{-c} |
\lambda \bigr]
\\[-2pt]
& \le &(\psimax\lambda_e + \lambda_{u_i})^c \Ex
\bigl[ \bigl( u_i^T \bigl(R_i
Q_{\lambda}^{-1}R_i^T
\bigr)^{-1} u_i \bigr)^{-c} | \lambda \bigr]
\\[-2pt]
& \le&2^{-c} \frac{\Gamma(\sfrac{q_i}{2} - c)}{\Gamma(\sfrac{q_{i}}{2})} (\psimax \lambda_e +
\lambda_{u_i})^c
\\[-2pt]
& \le&2^{-c} \frac{\Gamma(\sfrac{q_i}{2} - c)}{\Gamma(\sfrac{q_{i}}{2})} \bigl(\psimax^c
\lambda_e^c + \lambda_{u_i}^c\bigr).
\end{eqnarray*}
Define $\delta'(\cdot)$ as follows:
\[
\delta'(\alpha) = \frac{\psimax^c}{\alpha} \sum
_{i \in B} \frac{\Gamma(a_i + \sfrac{q_i}{2} + c)} {\Gamma(a_i + \sfrac{q_i}{2})} \frac{\Gamma(\sfrac{q_i}{2} - c)}{\Gamma(\sfrac{q_{i}}{2})}.
\]
Now we have
\begin{eqnarray*}
\sum_{i \in B} \Ex \bigl[ \tilde{\lambda}_{u_i}^c|
\lambda \bigr] & \le& \sum_{i \in B} \frac{\Gamma(a_i + \sfrac{q_i}{2} + c)} {\Gamma(a_i +
\sfrac{q_i}{2})}
\frac{\Gamma(\sfrac{q_i}{2} -
c)}{\Gamma(\sfrac{q_{i}}{2})} \bigl(\psimax^c \lambda_e^c
+ \lambda_{u_i}^c\bigr)
\\
& = &\alpha\delta'(\alpha) \lambda_e^c +
\sum_{i
\in B} \frac{\Gamma(a_i + \sfrac{q_i}{2} + c)} {\Gamma(a_i +
\sfrac{q_i}{2})} \frac{\Gamma(\sfrac{q_i}{2} -
c)}{\Gamma(\sfrac{q_{i}}{2})}
\lambda_{u_i}^c.
\end{eqnarray*}
Next, defining
\[
\rho'(\alpha) = \max \biggl\{ \delta'(\alpha), \max
_{i \in B} \biggl\{ \frac{\Gamma(a_i + \sfrac{q_i}{2} + c)} {\Gamma(a_i + \sfrac{q_i}{2})} \frac{\Gamma(\sfrac{q_i}{2} - c)}{\Gamma(\sfrac{q_{i}}{2})} \biggr\}
\biggr\},
\]
we have
\[
\sum_{i \in B} \Ex \bigl[ \tilde{\lambda}_{u_i}^c|
\lambda \bigr] \le \rho'(\alpha) \biggl( \alpha\lambda_e^c
+ \sum_{i \in B} \lambda_{u_i}^c
\biggr).
\]
Hence, all we have left to do is to prove that there exist $c \in
(0,1/2)$ and $\alpha$ satisfying \eqref{eq:alpha_bound} such that
$\rho'(\alpha) \in[0,1)$. First, define $\tilde{a} = -\max_{i \in B}
a_i$, and note that this quantity is positive. R\&H12 show that, if
$c \in(0,1/2) \cap(0,\tilde{a})$, then
\[
\max_{i \in B} \biggl\{ \frac{\Gamma(a_i + \sfrac{q_i}{2} + c)}{
\Gamma(a_i + \sfrac{q_i}{2})} \frac{\Gamma(\sfrac{q_i}{2} -
c)}{\Gamma(\sfrac{q_{i}}{2})}
\biggr\} < 1.
\]
Fix $c \in(0,1/2) \cap(0,\tilde{a})$. Now it suffices to show that
there exists an $\alpha$ satisfying \eqref{eq:alpha_bound} such that
$\delta'(\alpha) < 1$. But $\delta'(\alpha) < 1$ as long as
\[
\alpha> \psimax^c \sum_{i \in B}
\frac{\Gamma(a_i + \sfrac{q_i}{2} +
c)} {\Gamma(a_i + s\sfrac{q_i}{2})} \frac{\Gamma(\sfrac{q_i}{2} -
c)}{\Gamma(\sfrac{q_{i}}{2})}.
\]
So, \eqref{eq:suffice} is satisfied for $c \in(0,1/2) \cap
(0,\tilde{a})$ and
\begin{eqnarray*}
\alpha &>& \max \biggl\{ \frac{\sum_{i=1}^{r} \frac{\Gamma(a_{i} +
\sfrac{q_{i}}{2} - s)}{\Gamma( a_{i} + \sfrac{q_{i}}{2})}  (
\sfrac{\trace(R_{i}(I - P_{Z^{T}Z})R_{i}^{T})}{2}
 )^s}{1-\frac{\Gamma(a_{0} + \sfrac{N}{2} -s)}{\Gamma(a_{0} +
\sfrac{N}{2})}  ( \sfrac{\rank(Z)}{2}  )^s},
\\
&&{}  \psimax^c \sum
_{i
\in B} \frac{\Gamma(a_i + \sfrac{q_i}{2} + c)} {\Gamma(a_i +
\sfrac{q_i}{2})} \frac{\Gamma(\sfrac{q_i}{2} -
c)}{\Gamma(\sfrac{q_{i}}{2})} \biggr\}.
\end{eqnarray*}\upqed
\end{pf*}

\section{Discussion}
\label{sec:dis}

Our Corollary \ref{cor1} is a direct generalization of Rom{\'a}n and Hobert's \cite{romahobe2015}
Proposition 1 where we have removed all restrictions on the matrix
$X$. We now present a related result from \cite{abra2015} that is
established using a different drift function.

%
\begin{proposition}
\label{prop:alt}
Under a proper prior, the block Gibbs Markov chain,
$\{(\lambda_n,\theta_n)\}_{n=0}^\infty$, is geometrically ergodic if
$\min\{a_0, a_1, \dots, a_r \}>1$.
\end{proposition}

Like Corollary \ref{cor1}, this result holds for any $X$. Neither result is
uniformly better than the other. That is, there are situations where
the conditions of Corollary \ref{cor1} hold, but those of
Proposition~\ref{prop:alt} do not, and vice versa. However, the
condition $\min\{a_0, a_1, \dots, a_r \}>1$ appears to be more
restrictive than the conditions of Corollary \ref{cor1} in nearly all
\textit{practical} settings. In fact, the only examples we could find
where Proposition~\ref{prop:alt} is better than Corollary \ref{cor1} involve
models that have more random effects than observations. On the other
hand, we do feel that Proposition~\ref{prop:alt} is worth mentioning
because its simple form may render it useful to practitioners. For
example, in an exploratory phase where a number of different models
are being considered for a given set of data, one could avoid having
to recheck the conditions of Corollary \ref{cor1} each time the model is
changed simply by taking $a_0 = a_1 = \cdots= a_r = a > 1$ for all
models under consideration.

\begin{appendix}\label{append}

\section{Preliminary results}\label{app:prel}

Let $k = \rank(\tilde{X}) = \operatorname{rank}(X) \leq\min\{N,p\}$, and
consider a singular value decomposition of $\tilde{X}$ given by $U D
V^{T}$, where $U$ and $V$ are orthogonal matrices of dimension $N$ and
$p$, respectively, and
\[
D:= \lleft[ %
\begin{array} {cc} D_{*} & 0_{k, p-k}
\\
0_{N-k, k} & 0_{N-k, p-k} \end{array} %
 \rright],
\]
where $D_{*}:= \operatorname{diag}\{d_{1},\ldots,d_{k}\}$. The values
$d_{1},\ldots,d_{k}$ are the singular values of $\tilde{X}$, which are
strictly positive. Again, $\dmax$ denotes the largest singular value.
The following result is an extension of Lemmas 4 and 5 in R\&H15.

%
\begin{lemma}
\label{lem:rep1}
The matrix $\msl$ can be represented as $U \hsl U^T$ where
$H_{\lambda}$ is an $N\times N$ diagonal matrix, $\hsl=
\operatorname{diag}\{h_1,\ldots,h_N\}$, where
\[
h_i = \cases{ \displaystyle\frac{1}{\lambda_e d_i^2+1}, &\quad $i \in\{1,\ldots,k\}$,
\cr
1, &\quad $i \in\{k+1,\ldots, N\}$.}
\]
Furthermore, $(\lambda_e \dmax^2 +1)^{-1} I \preceq\msl\preceq I$.
\end{lemma}

\begin{pf}
Using the definitions of $T_{\lambda}^{-1}$ and $\tilde{X}$, we have
\[
\msl= I -\lambda_e X T_{\lambda}^{-1}X^T
= I - \lambda_e \tilde{X}\bigl(\lambda_e
\tilde{X}^T \tilde{X} + I\bigr)^{-1} \tilde{X}^T.
\]
Now using $\tilde{X} = U D V^{T}$ leads to
\[
\msl= U \bigl( I - \lambda_e D \bigl(\lambda_e
D^T D + I\bigr)^{-1} D^T \bigr) U^T.
\]
The matrix $\lambda_e D (\lambda_e D^T D + I)^{-1} D^T$ is an $N
\times N$ diagonal matrix whose $j$th diagonal element is given by
\[
\frac{\lambda_e d_j^2}{\lambda_e d_j^2+1} I_{\{1,2,\dots,k\}}(j).
\]
Hence, $I - \lambda_e D (\lambda_e D^T D + I)^{-1} D^T = \hsl$, and
$\msl= U \hsl U^T$. To prove the second part, note that, for
$j=1,\dots,N$, $0 < (\lambda_e \dmax^2+1)^{-1} \le h_i \le1$. Thus,
\[
\bigl(\lambda_e \dmax^2+1\bigr)^{-1} I = U
\bigl(\lambda_e \dmax^2+1\bigr)^{-1}
U^T \preceq U \hsl U^T \preceq U U^T = I.
\]\upqed
\end{pf}

Next, we develop an extension of Lemma 2 in R\&H15. Define $\tilde{Z}
= U^T Z$, $\tilde{y} = U^T y$ and $\eta= V^T \Sigma_{\beta}^{-1/2}
\mu_{\beta}$. Also, let $\tilde{z}_i$ denote the $i$th column of
$\tilde{Z}^T$, and let $\tilde{y}_i$ and $\eta_i$ represent the $i$th
components of the vectors $\tilde{y}$ and $\eta$, respectively. Let
$t_1, t_2, \dots, t_{N+q}$ be a set of $q$-vectors defined as follows.
For $j=1,\dots,N$, let $t_j = \tilde{z}_j$, and let $t_{N+1},\dots,
t_{N+q}$ be the standard basis vectors in~$\mathbb{R}^q$. For
$i=1,\dots,N$, define
\[
C_i^* = \Biggl[ \sup_{a \in\mathbb{R}_{+}^{N+q}} t_i^T
\Biggl( t_i t_i^T + \sum
_{j \in\{1,2,\dots,N\} \setminus\{i\}} a_j t_j t_j^T
+ \sum_{j=N+1}^{N+q} a_j
t_j t_j^T + a_i I
\Biggr)^{-2} t_i \Biggr]^{\sfrac{1}{2}}.
\]
The $C_i^*$s are finite by \cite{kharhobe2011}, Lemma 3.

%
\begin{lemma}
\label{lemma:bounding_ineq}
For all $\lambda\in\mathbb{R}^{r+1}$,
\[
\bigl\llVert \lambda_e Q_{\lambda}^{-1}Z^T
\msl y\bigr\rrVert \leq \sum_{j=1}^N |
\tilde{y}_j| C_j^* < \infty
\]
and
\[
\bigl
\llVert \lambda_e Q_{\lambda
}^{-1}Z^T X
T_{\lambda}^{-1} \Sigma_\beta^{-1}
\mu_{\beta}\bigr\rrVert \leq\sum_{j=1}^k
d_j |\eta_j| C_j^* < \infty.
\]
\end{lemma}

\begin{pf}
Even though R\&H15 assume $X$ to be full column rank, their argument
still works to establish the first inequality, so we omit this
argument. We now establish the second inequality. First,
\[
U^T X T_{\lambda}^{-1}= U^T \tilde{X}
\bigl(\lambda_e \tilde{X}^T \tilde {X} + I
\bigr)^{-1} \Sigma_{\beta}^{1/2} = D\bigl(
\lambda_e D^T D +I\bigr)^{-1} V^T
\Sigma_{\beta}^{1/2}.
\]
Define $R_\lambda= D(\lambda_e D^T D +I)^{-1}$. This is an $N \times
p$ diagonal matrix, with diagonal elements
$r_1,r_2,\dots,r_{\min\{N,p\}}$. These take the form
\[
r_j = \frac{d_j}{\lambda_e d_j^2 + 1} I_{\{1,2,\dots,k\}}(j).
\]
Now
\begin{eqnarray*}
\bigl\llVert \lambda_e Q_{\lambda}^{-1}Z^T
X T_{\lambda}^{-1}\Sigma _\beta ^{-1}
\mu_{\beta}\bigr\rrVert & =& \bigl\llVert \lambda_e
Q_{\lambda}^{-1}\tilde{Z}^T R_\lambda
V^T \Sigma _\beta^{-1/2}\mu_{\beta}\bigr
\rrVert
\\
& =& \bigl\llVert \lambda_e \bigl( \lambda_e
Z^T \msl Z + \Lambda \bigr)^{-1} \tilde{Z}^T
R_\lambda V^T \Sigma_\beta^{-1/2}\mu
_{\beta }\bigr\rrVert
\\
& =& \bigl\llVert \bigl( \tilde{Z}^T \hsl\tilde{Z} +
\lambda_e^{-1} \Lambda \bigr)^{-1}
\tilde{Z}^T R_\lambda\eta\bigr\rrVert
\\
& =& \Biggl\llVert \sum_{i=1}^k \bigl(
\tilde{Z}^T \hsl\tilde{Z} + \lambda_e^{-1}
\Lambda \bigr)^{-1} \tilde{z}_i r_i
\eta_i \Biggr\rrVert
\\
& \le &\sum_{i=1}^k \bigl\llVert \bigl(
\tilde{Z}^T \hsl\tilde{Z} + \lambda_e^{-1}
\Lambda \bigr)^{-1} \tilde{z}_i r_i
\eta_i \bigr\rrVert
\\
& = &\sum_{i=1}^k \Biggl\llVert \Biggl(
\sum_{j=1}^N \tilde{z}_j
\tilde{z}_j^T h_j + \lambda_e^{-1}
\Lambda \Biggr)^{-1} \tilde{z}_i r_i
\eta_i \Biggr\rrVert
\\
& = &\sum_{i=1}^k \biggl\llVert \biggl(
\tilde{z}_i \tilde{z}_i^T + \sum
_{j
\neq i} \tilde{z}_j \tilde{z}_j^T
\frac{h_j}{h_i} + h_i^{-1} \lambda_e^{-1}
\Lambda \biggr)^{-1} \tilde{z}_i \frac{r_i}{h_i}
\eta_i \biggr\rrVert
\\
& = &\sum_{i=1}^k d_i |
\eta_i| \biggl\llVert \biggl( \tilde{z}_i
\tilde{z}_i^T + \sum_{j \neq i}
\tilde{z}_j \tilde{z}_j^T \frac{h_j}{h_i}
+ h_i^{-1} \lambda_e^{-1} \Lambda
\biggr)^{-1} \tilde{z}_i \biggr\rrVert,
\end{eqnarray*}
where, in the last step, we have used the fact that $h_i d_i=r_i$ for
$i=1,\dots,k$. For $i=1,2,\dots,k$, define
\[
C_i(\lambda) = \biggl\llVert \biggl( \tilde{z}_i
\tilde{z}_i^T + \sum_{j
\neq i}
\tilde{z}_j \tilde{z}_j^T \frac{h_j}{h_i}
+ h_i^{-1} \lambda_e^{-1} \Lambda
\biggr)^{-1} \tilde{z}_i \biggr\rrVert.
\]
Define $\lambda_{\bullet} = \sum_{i=1}^r \lambda_{u_i}^{-1}$. Fix
$i$, and note that
\begin{eqnarray*}
C^2_i(\lambda) & =& \tilde{z}_i^T
\biggl( \tilde{z}_i \tilde{z}_i^T + \sum
_{j \neq i} \tilde{z}_j \tilde{z}_j^T
\frac{h_j}{h_i} + h_i^{-1} \lambda_e^{-1}
\Lambda \biggr)^{-2} \tilde{z}_i
\\
& = &\tilde{z}_i^T \biggl( \tilde{z}_i
\tilde{z}_i^T + \sum_{j \neq i}
\tilde{z}_j \tilde{z}_j^T \frac{h_j}{h_i}
+ h_i^{-1} \lambda_e^{-1} \bigl(
\Lambda- \lambda_{\bullet}^{-1} I\bigr) + \frac{1}{h_i \lambda_e
\lambda_{\bullet}} I
\biggr)^{-2} \tilde{z}_i.
\end{eqnarray*}
Define $\{w_j\}_{j=1}^{N+q}$ as follows:
%
%
\begin{equation}
w_j = \cases{ \displaystyle\frac{h_j}{h_i}, &\quad $j = 1,
\dots,i-1,i+1,\dots,N$,
\vspace*{3pt}\cr
\displaystyle\frac{1}{h_i \lambda_e \lambda_{\bullet}}, &\quad $j=i$,
\vspace*{3pt}\cr
\displaystyle\frac{\lambda_{u_1}-\lambda_{\bullet}^{-1}}{h_i \lambda_e}, &\quad $j = N+1,\dots,N+q_1$,
\vspace*{3pt}\cr
\displaystyle\frac{\lambda_{u_2}-\lambda_{\bullet}^{-1}}{h_i \lambda_e}, &\quad $j = N+q_1+1,
\dots,N+q_1+q_2$,
\vspace*{3pt}\cr
\vdots&\quad \vdots
\vspace*{3pt}\cr
\displaystyle
\frac{\lambda_{u_r}-\lambda_{\bullet}^{-1}}{h_i \lambda_e}, &\quad $j = N+q_1+\cdots+q_{r-1}+1,
\dots,N+q$.}
\end{equation}
Then
\[
C^2_i(\lambda) = t_i^T \Biggl(
t_i t_i^T + \sum
_{j \in\{1,2,\dots,N\}
\setminus\{i\}} w_j t_j t_j^T
+ \sum_{j=N+1}^{N+q} w_j
t_j t_j^T + w_i I
\Biggr)^{-2} t_i.
\]
Clearly, $w_j>0$ for all $j=1,\dots,N+q$. It follows that
\[
C^2_i(\lambda) \le\sup_{a \in\mathbb{R}_{+}^{N+q}}
t_i^T \Biggl( t_i t_i^T
+ \sum_{j \in\{1,2,\dots,N\} \setminus\{i\}} a_j t_j
t_j^T + \sum_{j=N+1}^{N+q}
a_j t_j t_j^T +
a_i I \Biggr)^{-2} t_i = \bigl(
C_i^* \bigr)^2.
\]
Hence,
\[
\bigl\llVert \lambda_e Q_{\lambda}^{-1}Z^T
X T_{\lambda}^{-1}\Sigma _\beta ^{-1}
\mu_{\beta}\bigr\rrVert \le\sum_{i=1}^k
d_i |\eta_i| C_i^*.
\]\upqed
\end{pf}

\section{Proof of Lemma~\texorpdfstring{\protect\ref{lem:ineq1}}{1}}
\label{app:lemma1}

\setcounter{lemma}{0}
\begin{lemma}
For each $i \in\{1,2,\dots,r\}$, we have
\[
\trace\bigl(R_{i} Q_{\lambda}^{-1}R_{i}^{T}
\bigr) \leq\bigl(\dmax^2 + \lambda_{e}^{-1}\bigr)
\trace \bigl( R_{i} \bigl(Z^{T}Z\bigr)^+
R_{i}^{T} \bigr) + \trace \bigl( R_{i}(I -
P_{Z^{T}Z}) R_{i}^{T} \bigr) \sum
_{j=1}^{r} \lambda_{u_j}^{-1}.
\]
\end{lemma}

\begin{pf}
From Lemma~\ref{lem:rep1} we have
\[
\qsl= \lambda_e Z^T \msl Z + \Lambda\succeq
\frac{\lambda_e}{\lambda_e \dmax^{2} + 1} Z^T Z + \Lambda\succeq \frac{\lambda_e}{\lambda_e \dmax^{2} + 1}
Z^T Z + \lmin I,
\]
where $\lmin= \min\{\lambda_{u_1},\dots,\lambda_{u_r}\}$. Letting
$O \Psi O^T$ be the spectral decomposition of $Z^T Z$, we have
%
%
\begin{equation}
\label{eq:Qineq} Q_{\lambda}^{-1}\preceq \biggl( \frac{\lambda_e}{\lambda_e \dmax
^{2} + 1}
Z^T Z + \lmin I \biggr)^{-1} = O \biggl( \frac{1}{\dmax^{2} +
\lambda_e^{-1}}
\Psi+ \lmin I \biggr)^{-1} O^T.
\end{equation}
Next, let $\Psi^+$ be a $q \times q$ diagonal matrix whose $i$th
diagonal element is
\[
\psi^+_i = \psi^{-1}_i \bigl(
1-I_{\{0\}}(\psi_i) \bigr).
\]
Now note that, for $i=1,\dots,q$, we have
\[
\biggl( \frac{\psi_i}{\dmax^{2} + \lambda_e^{-1}} + \lmin \biggr)^{-1} \le \bigl(
\dmax^{2} + \lambda_e^{-1} \bigr)
\psi_i^+ + \lmin^{-1} I_{\{0\}}\bigl(
\psi^+_i\bigr).
\]
Hence,
%
%
\begin{equation}
\label{eq:Iineq} \biggl( \frac{1}{\dmax^{2} + \lambda_e^{-1}} \Psi+ \lmin I \biggr)^{-1}
\preceq \bigl( \dmax^{2} + \lambda_e^{-1} \bigr)
\Psi^+ + \lmin^{-1} ( I - P_\Psi ),
\end{equation}
where $P_\Psi$ is a $q \times q$ diagonal matrix whose $i$th diagonal
entry is $1-I_{\{0\}}(\psi_i)$. Combining \eqref{eq:Qineq} and
\eqref{eq:Iineq} yields
\begin{eqnarray*}
Q_{\lambda}^{-1}&\preceq \bigl( \dmax^{2} +
\lambda_e^{-1} \bigr) O \Psi^+ O^T +
\lmin^{-1} O ( I - P_\Psi ) O^T
\\
&= \bigl( \dmax^{2} + \lambda_e^{-1} \bigr)
\bigl( Z^T Z\bigr)^+ + \lmin^{-1} O ( I - P_\Psi
) O^T.
\end{eqnarray*}
Let $\mathcal{I} = \{i \in\{1,\dots,q\}: \psi_i>0 \}$, and let
$\tilde{O}$ be the sub-matrix of $O$ consisting of the column vectors
$o_i$ where $i \in\mathcal{I}$. Then
\[
O P_{\Psi} O^T = \sum_{i \in\mathcal{I}}
o_i o_i^T = \tilde{O} \tilde{O}^T.
\]
Since $\{o_i\}_{i \in\mathcal{I}}$ forms an orthonormal basis for the
column space of $Z^T Z$, it follows that $\tilde{O} \tilde{O}^T$ is
the orthogonal projection onto $Z^T Z$. Consequently,
\[
O ( I-P_{\Psi} )O^T = O O^T-O P_{\Psi}
O^T = I - \tilde{O} \tilde{O}^T = I - P_{Z^T Z}.
\]
Thus,
\[
Q_{\lambda}^{-1}\preceq \bigl( \dmax^{2} +
\lambda_e^{-1} \bigr) \bigl( Z^T Z\bigr)^+ + (
I - P_{Z^T Z} ) \sum_{i=1}^r
\lambda_{u_i},
\]
and finally,
\[
\trace\bigl(R_{i} Q_{\lambda}^{-1}R_{i}^{T}
\bigr) \leq\bigl(\dmax^2 + \lambda_{e}^{-1}\bigr)
\trace \bigl( R_{i} \bigl(Z^{T}Z\bigr)^+
R_{i}^{T} \bigr) + \trace \bigl( R_{i}(I -
P_{Z^{T}Z}) R_{i}^{T} \bigr) \sum
_{j=1}^{r} \lambda_{u_j}^{-1}.
\]\upqed
\end{pf}

\section{Proof of Lemma~\texorpdfstring{\protect\ref{lem:ineq2}}{2}}
\label{app:lemma2}

\begin{lemma}
$\trace( W \Var(\theta|\lambda) W^{T})
\leq\lambda_{e}^{-1} \rank(Z) + \dmax^{2} \rank(Z) + \trace(X
\Sigma_{\beta} X^T)$.
\end{lemma}

\begin{pf}
R\&H15 show that
\[
\trace\bigl( W \Var(\theta|\lambda) W^T\bigr) = \trace\bigl( Z
Q_{\lambda
}^{-1}Z^T\bigr) + \trace\bigl( X
T_\lambda^{-1} X^T\bigr) - \trace \bigl( (I-\msl) Z
Q_{\lambda}^{-1}Z^T (I+\msl) \bigr),
\]
and that $\trace ( (I-\msl) Z Q_{\lambda}^{-1}Z^T (I+\msl)
) \ge0$.
Hence,
\[
\trace\bigl( W \Var(\theta|\lambda) W^T\bigr) \le\trace\bigl( Z
Q_{\lambda
}^{-1}Z^T\bigr) + \trace\bigl( X
T_\lambda^{-1} X^T\bigr),
\]
Next, note that $\Sigma_\beta^{-1}\preceq\lambda_e X^T X + \Sigma
_\beta^{-1}= \tsl$.
Hence, $\Sigma_{\beta} \succeq T_{\lambda}^{-1}$, and
\[
\trace\bigl( X T_\lambda^{-1} X^T\bigr) \le\trace
\bigl(X \Sigma_{\beta} X^T\bigr).
\]
Now, from Lemma~\ref{lem:rep1}, we have
\[
\frac{\lambda_e}{\lambda_e \dmax^2+1} Z^T Z + \Lambda\preceq \lambda_e
Z^T \msl Z + \Lambda= \qsl,
\]
and it follows that
\[
\trace\bigl(Z Q_{\lambda}^{-1}Z^T\bigr) \leq\trace
\biggl( Z \biggl( \frac{\lambda_e}{\lambda_e \dmax^2+1} Z^TZ + \Lambda
\biggr)^{-1} Z^T \biggr).
\]
Finally, using Lemma 3 from R\&H15, we have
\begin{eqnarray*}
\trace \biggl( Z \biggl( \frac{\lambda_e}{\lambda_e \dmax^2+1} Z^TZ + \Lambda
\biggr)^{-1} Z^T \biggr) & \leq &\biggl( \frac{\lambda_e}{\lambda_e \dmax^2+1}
\biggr)^{-1} \rank(Z)
\\
& =& \lambda_e^{-1} \rank(Z) + \dmax^2
\rank(Z).
\end{eqnarray*}\upqed
\end{pf}

\section{Proof of Lemma~\texorpdfstring{\protect\ref{lem:ineq3}}{3}}
\label{app:lemma3}

\begin{lemma}There exist finite constants $K_1$ and $K_2$,
not depending on $\lambda$, such that $\llVert \Ex[R_i u |
\lambda]\rrVert
\le\sqrt{q_i} K_1 $ for $i=1,\dots,r$, and $\llVert  y - W \Ex
[\theta| \lambda]\rrVert  \le K_2$.
\end{lemma}

\begin{pf}
From \eqref{eq:postmean} and Lemma~\ref{lemma:bounding_ineq}, we have
\begin{eqnarray*}
\bigl\llVert \Ex[u | \lambda]\bigr\rrVert & = &\bigl\llVert \lambda_e
Q_{\lambda}^{-1}Z^T \bigl(\msl y - X
T_{\lambda}^{-1}\Sigma_\beta ^{-1}
\mu_{\beta}\bigr)\bigr\rrVert
\\
& \le &\bigl( \bigl\llVert \lambda_e Q_{\lambda}^{-1}
Z^T \msl y\bigr\rrVert + \bigl\llVert \lambda_e
Q_{\lambda}^{-1}Z^T X T_{\lambda }^{-1}
\Sigma _\beta^{-1} \mu_{\beta}\bigr\rrVert \bigr)
\\
& \le& \Biggl( \sum_{j=1}^N \llvert
y_j\rrvert C_j^* + \sum_{j=1}^k
d_j |\eta_j| C_j^* \Biggr):=
K_1.
\end{eqnarray*}
Now, for each $i \in\{1,\dots,q\}$, we have
\[
\bigl\llVert \Ex[R_iu | \lambda]\bigr\rrVert \le\llVert
R_i\rrVert K_1 = \sqrt{\trace
\bigl(R_i^T R_i\bigr)} K_1 =
\sqrt{q_i} K_1.
\]
This proves the first part. Now, it follows from page 10 of R\&H15
that
\[
\bigl\llVert y - W \Ex[\theta| \lambda]\bigr\rrVert \le\llVert \msl\rrVert
\llVert y\rrVert + \bigl\llVert X T_{\lambda}^{-1}
\Sigma_\beta^{-1}\mu_{\beta}\bigr\rrVert + \llVert \msl
\rrVert \llVert Z\rrVert \bigl\llVert \Ex [u|\lambda]\bigr\rrVert.
\]
Now, using Lemma~\ref{lem:rep1}, and the fact that $h_i \le1$, for
$i=1,\dots,N$, we have
\[
\llVert \msl\rrVert ^2 = \trace\bigl(\msl^T \msl\bigr) =
\sum_{j=1}^N h_i^2
\le N.
\]
Recall from the proof of Lemma~\ref{lemma:bounding_ineq} that $U^T X
T_{\lambda}^{-1}= R_\lambda V^T \Sigma_{\beta}^{1/2}$, and note that
\[
\llVert R_\lambda\rrVert ^2 = \trace\bigl(R_\lambda^T
R_\lambda \bigr) = \sum_{j=1}^k
r_i^2 \le k \dmax^2.
\]
Therefore,
\begin{eqnarray*}
\bigl\llVert X T_{\lambda}^{-1}\Sigma_\beta^{-1}
\mu_{\beta}\bigr\rrVert & =& \bigl\llVert U U^T X
T_{\lambda}^{-1}\Sigma_\beta^{-1} \mu
_{\beta}\bigr\rrVert
\\
& =& \bigl\llVert U R_\lambda V^T \Sigma _{\beta}^{-1/2}
\mu_{\beta}\bigr\rrVert
\\
& \le&\llVert U\rrVert \llVert R_\lambda\rrVert \bigl\llVert
V^T \Sigma_{\beta}^{-1/2} \mu_{\beta}\bigr
\rrVert
\\
& \le&\sqrt{N} \sqrt{k} \dmax\bigl\llVert V^T \Sigma_{\beta}^{-1/2}
\mu_{\beta}\bigr\rrVert.
\end{eqnarray*}
Putting all of this together, we have
\begin{eqnarray*}
\bigl\llVert y - W \Ex[\theta| \lambda]\bigr\rrVert & \le&\llVert \msl\rrVert
\llVert y\rrVert + \bigl\llVert X T_{\lambda}^{-1}
\Sigma_\beta^{-1}\mu_{\beta}\bigr\rrVert + \llVert \msl
\rrVert \llVert Z\rrVert \bigl\llVert \Ex[u|\lambda]\bigr\rrVert
\\
& \le&\sqrt{N} \llVert y\rrVert + \sqrt{N} \sqrt{k} \dmax\bigl\llVert
V^T \Sigma_{\beta}^{-1/2} \mu_{\beta}\bigr
\rrVert + \sqrt{N} \llVert Z\rrVert K_1.
\end{eqnarray*}\upqed
\end{pf}

\section{Proof of Lemma~\texorpdfstring{\protect\ref{lem:ineq4}}{4}}
\label{app:lemma4}

\begin{lemma}
For each $i \in\{1,2,\dots,r\}$, we have
$(\psimax\lambda_e + \lambda_{u_i})^{-1} I_{q_i} \preceq R_i
Q_{\lambda}^{-1}
R_i^T $.
\end{lemma}

\begin{pf}
Lemma~\ref{lem:rep1} implies that $Z^T \msl Z \preceq Z^T Z$. It
follows that
\[
\qsl= \lambda_e Z^T \msl Z + \Lambda\preceq
\lambda_e Z^T Z + \Lambda \preceq\lambda_e
\psimax I + \Lambda.
\]
Thus,
\[
(\lambda_e \psimax+ \lambda_{u_i})^{-1} I =
R_i ( \lambda_e \psimax I + \Lambda)^{-1}
R_i^T \preceq R_i Q_{\lambda}^{-1}R_i^T.
\]\upqed
\end{pf}
\end{appendix}


\section*{Acknowledgment}
The authors thank three anonymous reviewers for helpful comments and suggestions that
led to a substantially improved version of the paper.

Tavis Abrahamsen supported by NSF Grant 08-01544 (in the Quantitative Spatial Ecology, Evolution
and Environment Program at the University of Florida).
James P. Hobert supported by NSF Grant DMS-11-06395.


%

\printhistory
\end{document}